\font\cmrfootnote=cmr10 scaled 750
\font\cmrhalf=cmr10 scaled \magstephalf

\font\cmrscfootnote=cmr7 scaled 750
\font\cmrschalf=cmr7 scaled \magstephalf

\font\cmrscscfootnote=cmr5 scaled 750
\font\cmrscschalf=cmr5 scaled \magstephalf

\font\mitfootnote=cmmi10 scaled 750
\font\mithalf=cmmi10 scaled \magstephalf

\font\mitscfootnote=cmmi7 scaled 750
\font\mitschalf=cmmi7 scaled \magstephalf

\font\mitscscfootnote=cmmi5 scaled 750
\font\mitscschalf=cmmi5 scaled \magstephalf

\font\cmsyfootnote=cmsy10 scaled 750
\font\cmsyhalf=cmsy10 scaled \magstephalf

\font\cmsyscfootnote=cmsy7 scaled 750
\font\cmsyschalf=cmsy7 scaled \magstephalf

\font\cmsyscscfootnote=cmsy5 scaled 750
\font\cmsyscschalf=cmsy5 scaled \magstephalf

\font\cmexfootnote=cmex10 scaled 750
\font\cmexhalf=cmex10 scaled \magstephalf

\font\cmexscfootnote=cmex10 scaled 750
\font\cmexschalf=cmex10 scaled \magstephalf

\font\cmexscscfootnote=cmex7 scaled 750
\font\cmexscschalf=cmex7

\def\mathfootnote{\textfont0=\cmrfootnote \textfont1=\mitfootnote \textfont2=\cmsyfootnote \textfont3=\cmexfootnote
        \scriptfont0=\cmrscfootnote \scriptscriptfont0=\cmrscscfootnote \scriptfont1=\mitscfootnote \scriptscriptfont1=\mitscscfootnote
        \scriptfont2=\cmsyscfootnote \scriptscriptfont2=\cmsyscscfootnote \scriptfont3=\cmexscfootnote \scriptscriptfont3=\cmexscscfootnote}
\def\mathhalf{\textfont0=\cmrhalf \textfont1=\mithalf \textfont2=\cmsyhalf \textfont3=\cmexhalf
        \scriptfont0=\cmrschalf \scriptscriptfont0=\cmrscschalf \scriptfont1=\mitschalf \scriptscriptfont1=\mitscschalf
        \scriptfont2=\cmsyschalf \scriptscriptfont2=\cmsyscschalf \scriptfont3=\cmexschalf \scriptscriptfont3=\cmexscschalf}

\let\mf=\mathfootnote

\font\Bbb=msbm10

\def\outlin#1{\hbox{\Bbb #1}}

\font\call=cmsy10 

\def\cal{\call}

\font\small=cmr5                       
\font\notsosmall=cmr7
\font\cmreight=cmr8
 
\font\bfeight=cmbx8

\font\cmrhalf=cmr10 scaled \magstephalf
       
\font\two=cmr10 scaled \magstep2

\font\sansfoot=cmss8
\font\sans=cmss10

\font\twosans=cmss10 scaled \magstep2

\font\caps=cmcsc10

\let\mf=\mfootnote

\def\q{\quad}

\def\cut{\hfill\break}  \def\newline{\cut}
\def\h#1{\hbox{#1}}

\def\eqa{\eqalign} \def\eqano{\eqalignno}

 \def\eps{\epsilon}

 \def\L{\Lambda} 
\def\vp{\varphi}
  
\def\o{\omega}   \def\O{\Omega}

 \def\calH{{\h{\cal H}}}
\let\H=\calH

\def\calO{{\h{\cal O}}} 
 
\def\calL{{\h{\cal L}}} \def\L{{\calL}} 
\def\calI{{\h{\cal I}}}

 \def\RR{{\outlin R}} \def\ZZ{{\outlin Z}} \def\NN{{\outlin N}}
  \def\PP{{\outlin P}}

\chardef\dotlessi="10  
\chardef\inodot="10

\def\polishl{\char'40l}   
\def\polishL{\leavemode\setbox0=\hbox{L}\hbox to\wd0{\hss\char'40L}}

\def\diam{\hbox{\rm diam}}

\def\Aut{\h{\rm Aut}}  \def\aut{\h{\rm aut}}

   \def\Ric{\h{\rm Ric}}

\def\V{\frac1V}

\def\gij{g_{i\bar j}}

\def\AutMJ{\hbox{$\hbox{\rm Aut}(M,\hbox{\rm J})$}}

\def\opsh{$\o$-plurisubharmonic }

\def\GL2nR{\h{$GL(2n,\RR)$}}
\def\Sp2nR{\h{$Sp(2n,\RR)$}}

\def\dbz{d\bar z}

\def\dzidzjb{dz^i\w\dbz^j}

\def\part#1{\frac{\partial#1}{\partial t}}

\def\frac#1#2{{{#1}\over{#2}}}

\def\all{\forall}

\def\sm{\setminus}

\def\precpt
{\hbox{$\mskip3mu\mathhalf\subset\raise0.92pt\hbox{$\mskip-10mu\!\!\!\!
\mathfootnote\subset$}\mskip5mu$}}

\def\supsetnoteq{\hbox{$\mskip3mu\supset\raise-5.97pt
\hbox{$\mskip-10mu\!\!\!\scriptstyle\not=$}\mskip8mu$}}

\def\subsetnoteq{\hbox{$\mskip3mu\subset\raise-5.97pt
\hbox{$\mskip-11mu\!\!\scriptstyle\not=$}\mskip8mu$}}

\def\sseq{\subseteq}

\def\*{\star}

\def\w{\wedge}

\def\D{\Delta}
\def\N{\nabla} 
\def\dbar{\bar\partial}
\def\del{\partial} \def\pa{\partial}
\def\ddbar{\partial\dbar}

\def\intM{\int_M}  \def\intm{\int_M}
\def\V{\frac1V}
\def\VintM{\V\int_M}

\def\i{\sqrt{-1}}

\def\ra{\rightarrow}

\def\arrow#1{\hbox to #1pt{\rightarrowfill}}

\def\thhnotsosmall#1{${\hbox{#1}}^{\hbox{\small th}}$}

\def\MA{Monge-Amp\`ere }

\def\K{K\"ahler }

\def\KE{K\"ahler-Einstein }  
\def\KR{K\"ahler-Ricci }

\def\Kolodziej{Ko\polishl{}odziej}

\def\po1{partition of unity }

\def\Cinf{\h{\cal C}^\infty}

\def\Linf{L^\infty}

\def\Loneloc{L_{\h{\small loc}}^1}

\def\CinfM{\Cinf(M)}

\def\strutdepth{\dp\strutbox}
\def\specialstar{\vtop to \strutdepth{
    \baselineskip\strutdepth
    \vss\llap{$\star$\ \ \ \ \ \ \ \ \  }\null}}
\def\marginalstar{\strut\vadjust{\kern-\strutdepth\specialstar}}
\def\marginal#1{\strut\vadjust{\kern-\strutdepth
    {\vtop to \strutdepth{
    \baselineskip\strutdepth
    \vss\llap{{ \small #1 }}\null} 
    }}
    }

\newcount\subsectionitemnumber
\def\clearsubsectionitemnumber{\subsectionitemnumber=0\relax}

\newcount\subsubsubsectionnumber
\def\clearsubsubsubsectionnumber{\subsubsubsectionnumber=0\relax}
\def\subsubsubsection#1{
\bigskip\noindent
\global\advance\subsubsubsectionnumber by 1%
{\rm
 \the\subsectionnumber.\the\subsubsectionnumber.\the\subsubsubsectionnumber}
{
#1.}
}

\newcount\subsubsectionnumber
\def\clearsubsubsectionnumber{\subsubsectionnumber=0\relax}
\def\subsubsection#1{
\clearsubsubsubsectionnumber
\bigskip\noindent
\global\advance\subsubsectionnumber by 1%
{%
\it \the\subsectionnumber.\the\subsubsectionnumber}
{
\it #1.}
}
\newcount\subsectionnumber
\def\clearsubsectionnumber{\subsectionnumber=0\relax}
\def\subsection#1{
\clearsubsectionitemnumber
\clearsubsubsectionnumber
\medskip\medskip\smallskip\noindent \global\advance\subsectionnumber by 1%
{%
\bf \the\subsectionnumber} 
{
\bf #1.}
}
\newcount\sectionnumber
\def\clearsectionnumber{\sectionnumber=0\relax}
\def\section#1{
\clearsubsectionnumber
\bigskip\bigskip\noindent \global\advance\sectionnumber by 1%
{%
\two
\the\sectionnumber} 
{
\two #1.}
}

\clearsectionnumber
\clearsubsectionnumber
\clearsubsubsectionnumber
\clearsubsubsubsectionnumber
\clearsubsectionitemnumber

\newcount\itemnumber
\def\clearitemnumber{\itemnumber=0\relax}
\def\c#1{ {\noindent\bf \the\itemnumber.} p. $#1$ \global\advance\itemnumber by 1}
\def\cn{ {\noindent\bf \the\itemnumber.} \global\advance\itemnumber by 1}

\clearitemnumber

\def\subsectionno#1{
\medskip\medskip\smallskip\noindent%
{%
\bf #1.}
}

\def\last#1{\advance\eqncount by -#1(\the\eqncount)\advance\eqncount by #1}
\def\llast{\advance\eqncount by -1(\the\eqncount)\advance\eqncount by 1}
\def\lllast{\advance\eqncount by -2(\the\eqncount)\advance\eqncount by 2}
\def\llllast{\advance\eqncount by -3(\the\eqncount)\advance\eqncount by 3}

\newcount\notenumber
\def\clearnotenumber{\notenumber=0\relax}
\def\note#1{\advance\notenumber by 1\footnote{${}^{\the\notenumber}$}%
{\lineskip0pt\notsosmall #1}}
\def\notewithcomma#1{\advance\notenumber by
1\footnote{${}^{\the\notenumber}$}
  {\lineskip0pt\notsosmall #1}}
                                                                                               
\clearnotenumber

\def\putnumber{%
\global\advance\subsectionitemnumber by 1{\the\subsectionnumber}.{\the\subsectionitemnumber}}
\def\numbering{{{\the\subsectionnumber}.{\the\subsectionitemnumber}}}

\def\FThm#1{\bigskip
%\vbox
{\noindent%
{\bf {T}heorem }
{\bf \putnumber.} {\it #1}\bigskip 
}}%\vglue1in}\vglue-1in}

\def\FCor#1{\bigskip
%\vbox
{\noindent%
{\bf {C}orollary }
{\bf \putnumber.} {\it #1}\bigskip
}}%\vglue1in}\vglue-1in}
\def\FLem#1{\bigskip
%\vbox
{\noindent%
{\bf {L}emma }
{\bf \putnumber.} {\it #1}\bigskip
}}%\vglue1in}\vglue-1in}

%\vglue1in}\vglue-1in}
%\vglue1in}\vglue-1in}

%\vglue0cm}\vglue-0cm}

%\vglue1in}\vglue-1in}

%\vglue1in}\vglue-1in}

%\vglue1in}\vglue-1in}

\def\FProb#1{\bigskip
%\vbox
{\noindent%
{\bf {P}roblem }
{\bf \putnumber.} {\it #1}\bigskip
}}%\vglue1in}\vglue-1in}

%\vglue1in}\vglue-1in}

\def\FDef#1{\bigskip
%\vbox
{\noindent%
{\bf {D}efinition }
{\bf \putnumber.} {\it #1}\bigskip
}}%\vglue1in}\vglue-1in}

\def\Abstract#1{
{\narrower\bigskip\bigskip
\noindent {{\bf Abstract.\ \ } #1}

}
}

\def\ref#1{{\bf[}{\sans #1}{\bf]}}
  % \ref with small skip

\def\reffoot#1{{\bfeight[}{\sansfoot #1}{\bfeight]}}

\def\opcit{\underbar{\phantom{aaaaa}}}

\def\sm{\smallskip}

\def\boxit#1{\vbox{\hrule\hbox{\vrule\kern3pt\vbox{\kern3pt#1\vglue3pt
\kern3pt}\kern3pt\vrule}\hrule}}

\long\def\frame#1#2#3#4{\hbox{\vbox{\hrule height#1pt
 \hbox{\vrule width#1pt\kern #2pt
 \vbox{\kern #2pt
 \vbox{\hsize #3\noindent #4}
\kern#2pt}
 \kern#2pt\vrule width #1pt}
 \hrule height0pt depth#1pt}}}

\def\kou{\frame{.3}{.5}{1pt}{\phantom{a}}}

\def\help{\ifmmode\aftergroup\noindent\quad\else\quad\fi}
\def\helpp{\ifmmode\aftergroup\noindent\hfill\else\hfill\fi}
\def\done{\helpp\kou\medskip}

\def\Ho{{\calH_{\o}}}

\def\on{\o^n}
\def\onminus1{\o^{n-1}/(n-1)!}
\def\ovp{\o_{\vp}}

\def\vpt{\vp_{\!t}}

\def\ovpt{\o_{\vpt}}

\def\ovpn{\o_{\vp}^n}
\def\ovptn{\o_{\vpt}^n}

\def\fo{f_{\o}}
\def\fovpt{f_{\ovpt}}

\font\call=cmsy10   \font\Bbb=msbm10
\def\outlin#1{\hbox{\Bbb #1}}
\def\H{\hbox{\call H}}  \def\Ho{\hbox{\call H}_\omega}
 
\def\Cinf{{\call C}^\infty}  \def\RR{{\outlin R}}
\def\frac#1#2{{{#1}\over{#2}}}
\newcount \eqncount
\def \eqnno{\global \advance \eqncount by 1 \futurelet \nexttok \parsenexttok}
\def \eqn{\global \advance \eqncount by 1 \eqno\futurelet \nexttok \parsenexttok}
\def \eqnd{\global\advance \eqncount by 1 \futurelet\nexttok\parsenexttokd}
\def \parsenexttok{\ifx \nexttok $\Nomark\else\expandafter \Mark\fi}
\def \parsenexttokd{\ifx \nexttok \hfil\Nomark\else\expandafter \Mark\fi}
\def \Nomark {(\the \eqncount)}\def \Mark #1{\xdef #1{(\the \eqncount)}#1}

\let\exa\expandafter
\catcode`\@=11 
\def\CrossWord#1#2#3{%
\def\@x@{}\def\@y@{#2}
\ifx\@x@\@y@ \def\@z@{#3}\else \def\@z@{#2}\fi
\exa\edef\csname cw@#1\endcsname{\@z@}}
\openin15=\jobname.ref
\ifeof15 \immediate\write16{No file \jobname.ref}%
   \else \input \jobname.ref \fi 
\closein15
\newwrite\refout
\openout\refout=\jobname.ref 
\def\warning#1{\immediate\write16{1.\the\inputlineno -- warning --#1}}
\def\Ref#1{%
\exa \ifx \csname cw@#1\endcsname \relax
\warning{\string\Ref\string{\string#1\string}?}%
    \hbox{$???$}%
\else \csname cw@#1\endcsname \fi}
\def\Tag#1#2{\begingroup
\edef\head{\string\CrossWord{#1}{#2}}%
\def\writeref{\write\refout}%
\exa \exa \exa
\writeref \exa{\head{\the\pageno}}%
\endgroup}

\catcode`\@=12

\def\Tagg#1#2{\Tag{#1}{#2}}  

\def\TaggThm#1{\Tagg{#1}{Theorem~\numbering}}
\def\TaggLemm#1{\Tagg{#1}{Lemma~\numbering}}

\def\TaggCor#1{\Tagg{#1}{Corollary~\numbering}}

\def\TaggSection#1{\Tagg{#1}{Section~\the\subsectionnumber}}
\def\TaggS#1{\Tagg{#1}{\S~\the\subsectionnumber}}
\def\TaggSubS#1{\Tagg{#1}{\hbox{\S\S\unskip\the\subsectionnumber.\the\subsubsectionnumber}}}
\def\TaggSubsection#1{\Tagg{#1}{Subsection~\the\subsectionnumber.\the\subsubsectionnumber}}
\def\TaggEq#1{\Tagg{#1}{(\the\eqncount)}}
\def\Taggf#1{ \Tagg{#1}{{\bf[}{\sans #1}{\bf]}} }
\def\TaggDef#1{\Tagg{#1}{Definition~\numbering}}

\def\JJJ{\h{\rm J}}
\def\Ric{\hbox{\rm Ric}\,}

\def\Hc{\H_{c_1}}

\def\O{\Omega}

\def\Rico{\Ric\o}

\def\HO{\H_{\Omega}}
\def\Ho{{\calH_{\o}}}

\def\MA{Monge-Amp\`ere }
\def\KR{K\"ahler-Ricci }

\def\Vm{V^{-1}}  \def\V{\frac1V}

\def\fo{f_{\o}}

\magnification=1100
\hoffset1.1cm
\voffset1.3cm
\hsize5.09in
\vsize6.91in

\headline={\ifnum\pageno>1{\ifodd\pageno \oddheadline\else\evenheadline\fi}\fi}
\def\oddheadline{\centerline{\caps On Nadel multiplier ideal sheaves}}
\def\evenheadline{\centerline{\caps Y. A. Rubinstein}}

\overfullrule0pt
\parindent12pt
%\vglue0.6cm
%\vglue0.6cm

%\vglue-2cm
%{\sl To appear in Transactions of the American Mathematical Society}
%\vglue1.2cm

\noindent
\centerline{\twosans On the construction of Nadel multiplier ideal}
\medskip
\centerline{\twosans sheaves and the limiting behavior of the Ricci flow}
\smallskip

\font\tteight=cmtt8 

\bigskip
\centerline{\sans Yanir A. Rubinstein%
\footnote{$^*$}{\cmreight Massachusetts Institute of Technology. 
Email: {\tteight yanir@member.ams.org}
\hfill\break
\vglue-0.5cm
\hglue-\parindent\cmreight
Current address: Department of Mathematics, Princeton University, Princeton, NJ 08544.}}

\def\centereps#1#2#3{\vglue#2\relax\centerline{\hbox to#1%
            {\special{eps:#3.eps x=#1 y=#2}\hfil}}}

\font\cmreight=cmr8

%\vglue1cm

%\indent
%{\it To Aynat Rubinstein}

\vglue0.03in

\footnote{}{\hglue-\parindent\cmreight August \thhnotsosmall{12}, 2007. 
Revised October 2007. 
\hfill\break 
} %  
%       % 
    
\vglue-0.12in

\bigskip
\Abstract{
In this note we construct Nadel multiplier ideal sheaves using the Ricci flow
on Fano manifolds. This extends a result of Phong, \v Se\v sum and Sturm.
These sheaves, like their counterparts constructed by Nadel for the continuity method, can be used to obtain
an existence criterion for K\"ahler-Einstein metrics. 
}

\bigskip

\medskip
\subsectionnumber=0\relax

\subsection{Introduction}
\TaggSection{SectionIntroduction}
In this note we construct Nadel multiplier ideal sheaves on Fano
manifolds that do not admit \KE metrics, using the Ricci flow. 
The result is a simple consequence of
the uniformity of the Poincar\'e and Sobolev inequalities along the flow.
This allows to obtain another proof of the convergence of the Ricci
flow on a certain class of Fano manifolds.

The theory of obstructions to the existence of canonical \K metrics 
(see, e.g., \ref{F})
has a long history starting with the observation that a \KE manifold must have a definite or zero
first Chern class. Lichnerowicz and Matsushima  proved that for a constant
scalar curvature
\K manifold the group of automorphisms $\AutMJ$ is a complexification of the group of 
isometries. Later it was shown that on such a manifold the closed 1-form on the space of \K metrics
$\HO$  
(see the end of this section for notation and definitions)
defined by the scalar curvature minus its average must be basic with respect to $\AutMJ_0$, 
that is admit an $\AutMJ_0$-invariant potential function, or equivalently
its Futaki character must vanish. 
Kobayashi and L\"ubke proved that the tangent bundle of a \KE manifold is semistable.

In contrast to these necessary conditions, came work on certain sufficient conditions for the
existence of \KE metrics on Fano manifolds. 
Siu, Tian and Yau showed that certain finite groups of symmetries can be used to this end and 
produced several examples \ref{Si,T1,TY}. 
Tian studied the singular locus of sequences of plurisubharmonic functions
and introduced a sufficient condition in terms
of a related holomorphic invariant \ref{T1}: $\alpha_M(G)>n/(n+1)$, where $G$ is compact subgroup of $\AutMJ$ and
$$
\alpha_M(G)=\sup\,\{\;\alpha \,:\, \intm e^{-\alpha(\vp-\sup\vp)}\on<C_\alpha,\q \all\,\vp\in\Ho(G),\q[\o]=c_1\}.
$$
Then, Nadel introduced the notion of a multiplier
ideal sheaf on a compact \K manifold and showed that the nonexistence of certain such
$G$-invariant sheaves implies that $\alpha_M(G)\ge 1$ \ref{N}. 
This construction is related to a theory introduced earlier by Kohn in a different context \ref{K}. 

Tian translated the failure of the continuity method when 
a \KE metric does not
exist to the statement that a certain subsequence of \K potentials along Aubin's continuity path \ref{A2}
$$
\o_{\vp_t}^n=\on e^{f_\o-t\vp_t},\q t\in[0,t_0),\; t_0\le1, \eqn\AubinPathEq
$$ 
will diverge 
along a subvariety in a manner that can be ruled out when $\alpha_M(G)$ is large enough.
Nadel showed that furthermore the blow-up will occur along a subscheme cut out by a 
coherent sheaf of multiplier ideals satisfying certain cohomology vanishing 
conditions. These results provide a powerful tool in showing existence of \KE metrics, 
since these conditions are often violated in specific examples.
This technique was revisited by Demailly and Koll\'ar who also extended it to orbifolds \ref{DK}.

Nadel's main result can be stated as follows.%
\note{\cmreight
This theorem makes use of Nadel's vanishing theorem; see \Ref{MainSection} for the statement.
Also,  
Nadel's original statement includes (for simplicity) only the case $\mf\gamma\in(\frac n{n+1},1)$, however
later he comments on the possibility of extending the allowed interval for $\mf\gamma$ \reffoot{N, p. 579}.
}

\FThm{
\TaggThm{NadelThm}
\ref{N}
Let $(M,\JJJ)$ be a Fano manifold not admitting a \KE metric. 
Let $\gamma\in(\frac n{n+1},\infty)$ and let $\o\in\Hc$.
Then there exists a subsequence $\{\vp_{t_j}\}_{j\ge0}$ of solutions of \AubinPathEq\
such that $\vp_{t_j}-\sup\vp_{t_j}$ converges in the $L^1(M,\o)$-topology
to $\vp_\infty\in PSH(M,\JJJ,\o)$ and
$\calI(\gamma\vp_\infty)$ is a proper multiplier ideal sheaf
satisfying
$$
H^r(M,\calI(\gamma\vp_\infty)\otimes K_M^{-\lfloor\gamma\rfloor})=0, \q \all\, r\ge 1.
\eqn\NadelVanishingFanoEq
$$
}
\noindent
Following Nadel, the sheaves with $\gamma<1$ will be referred to as Nadel sheaves 
(see \Ref{NadelSheafDef} below).

The Ricci flow, introduced by Hamilton \ref{H1}, provides another method 
for constructing \KE metrics on a Fano manifold, and it is therefore 
natural to ask whether this method will also yield multiplier ideal obstruction
sheaves in the absence of a \KE metric.
It may be written as a flow equation on the space of \K potentials $\Ho$,
$$
\ovptn=\on e^{\fo-\vpt+\dot\vpt},\q\,\vp(0)=\vp_0.\eqn\FlowPotentialMAEq
$$
This question was first addressed by Phong, \v Se\v sum and Sturm
who proved the following result. 

\vglue-0.15cm
\FThm{
\TaggThm{PSSThm}
\ref{PSS}
Let $(M,\JJJ)$ be a Fano manifold not admitting a \KE metric. 
Let $\gamma\in(1,\infty)$ and let $\o\in\Hc$.
Then there exists an initial condition $\vp_0$ and a subsequence 
$\{\vp_{t_j}\}_{j\ge0}$ of solutions of 
\FlowPotentialMAEq\
such that $\vp_{t_j}$ converges in the $L^1(M,\o)$-topology to $\vp_\infty\in PSH(M,\JJJ,\o)$ and
$\calI(\gamma\vp_\infty)$ is a proper multiplier ideal sheaf satisfying 
\NadelVanishingFanoEq.
}

\vglue-0.15cm

The sheaves thus obtained require the exponent to lie in a smaller interval than in 
\Ref{NadelThm}, that is, they are not Nadel sheaves. This is a crucial difference between these results.
Indeed, the smaller the exponent $\gamma$, the stronger the vanishing theorem satisfied
by the corresponding multiplier ideal sheaf. 
The gist of \Ref{NadelThm} is to use functions invariant under a compact (not necessarily maximally compact) 
subgroup
of automorphisms $G$ in order to obtain $G$-invariant
subschemes satisfying in addition certain cohomology vanishing restrictions. When
$\gamma\in(n/(n+1),1)$, the sheaves constructed by Nadel are not only $G$-invariant, but also
satisfy 
$$
H^r(M,\calI(\gamma\vp_\infty))=0, \q \all\, r\ge 0,\eqn\UsualNadelVanishingEq
$$
while this need not hold for 
sheaves corresponding to exponents $\gamma>1$. The subschemes cut out by sheaves
satisfying \UsualNadelVanishingEq\ satisfy various restrictions, and
it is these restrictions that render the continuity method useful in proving 
existence of \KE metrics on
a large class of manifolds. For example, to state the simplest restrictions, 
Nadel shows that \UsualNadelVanishingEq\ implies
that the corresponding subscheme is connected and has arithmetic genus zero 
and that if it is one-dimensional it
is a tree of rational curves.
It is not clear how to use the sheaves with $\gamma>1$ to prove such results.

The main result of this note is that the Ricci flow does produce Nadel sheaves, 
with $\gamma\in(n/(n+1),1)$.
\vglue-0.1cm
\FThm{
\TaggThm{MainThm}
Let $(M,\JJJ)$ be a Fano manifold not admitting a \KE metric. 
Let $\gamma\in(n/(n+1),\infty)$ and let $\o\in\Hc$.
Then there exists an initial condition $\vp_0$ and a subsequence 
$\{\vp_{t_j}\}_{j\ge0}$ of solutions of 
\FlowPotentialMAEq\
such that $\vp_{t_j}-\V\intm\vp_{t_j}$ converges in the $L^1(M,\o)$-topology to $\vp_\infty\in PSH(M,\JJJ,\o)$ and
$\calI(\gamma\vp_\infty)$ is a proper multiplier ideal sheaf satisfying 
\NadelVanishingFanoEq.
}

Perelman proved that when a \KE metric exists the Ricci flow will converge to it in the sense
of Cheeger-Gromov. Therefore the following corollary, due to Nadel and Perelman, is known:

\FCor{\ref{N,TZ}
\TaggCor{ConvergenceCor}
Let $(M,\JJJ)$ be a Fano manifold and let $G$ be a compact subgroup
of \AutMJ. Assume 
that $(M,\JJJ)$ does not admit a $G$-invariant Nadel sheaf
as in \Ref{NadelSheafDef}. Then  
the Ricci flow will converge in the sense of Cheeger-Gromov to 
a \KE metric starting with any $G$-invariant \K metric.
}
\vglue-0.1cm
Observe that \Ref{MainThm} allows to obtain another proof of this corollary. This proof 
differs from the one obtained by the combination of the results of Nadel and Perelman only
in the method of obtaining the $\Linf$ estimate (and not in the method of obtaining 
higher order estimates and convergence \ref{CT,PSS,TZ}).
It has the virtue of simultaneously proving
that a \KE metric exists and that the flow will converge to it (instead of first using the continuity
method to prove that a \KE metric exists and then using this fact via properness of 
an energy functional and \Kolodziej's results \ref{Ko} to obtain the $\Linf$ estimate). 

In particular, this gives another proof of the theorem that the two-sphere $S^2$ 
(and more generally, complex projective space $\PP^n$) may be 
uniformized using the Ricci flow\note{%
\cmreight Indeed, observe that the group of automorphisms generated by the rotations $\mf z\mapsto ze^{\i\theta}$ 
and the
inversion $\mf z\mapsto 1/z$ acts without fixed points on $\mf S^2$ while a 0-dimensional Nadel subscheme is a single reduced point.}
 (\ref{CLT,CT,Ch1,H2} and references therein). Still more generally, \Ref{ConvergenceCor} applies
also to symmetric toric Fano manifolds (see \Ref{So}).

We emphasize that no new techniques are needed
beyond those in the continuity method setting described above and our purpose in this note is
solely to point out this similarity between the limiting behavior of the Ricci flow and
that of the more classical continuity method.
We believe \Ref{MainThm} adds important information regarding the limiting behavior 
of the Ricci flow beyond that in \Ref{PSSThm}. 

The proof of \Ref{MainThm} differs from that of \Ref{PSSThm} in that
it follows the lines of the original continuity method results  \ref{A2,N,Si,T1} instead of appealing to  
results of \Kolodziej. This is also the key to obtaining the result for singularity exponents
in the whole interval $(n/(n+1),\infty)$. 
The crucial ingredient that makes this possible is that the relevant estimates
on the Ricci flow established by Perelman, Ye and Zhang, some of which appeared after the work
of Phong-\v Se\v sum-Sturm, allow one to adapt the continuity
method arguments to the setting of the flow.\note{\cmreight With the exception of the case $\mf n=1$ ($\mf M=S^2$) 
that
does not make use of Ye and Zhang's estimate.}  This is described in \Ref{MainSection}. 
\Ref{RemarksSection} contains some remarks, including a brief conjectural discussion on a possible extension of the 
result to the setting of constant scalar curvature metrics.

\vglue-0.21cm
\subsectionno{Setup and notation}
Let $(M,\JJJ)$ be a connected compact closed \K manifold of complex dimension $n$
and let $\O\in H^2(M,\RR)$ be a \K class.
Let $d=\del+\dbar$ and define the Laplacian $\D=-\dbar\circ\dbar^\star-\dbar^\star\circ\dbar$ with
respect to 
a Riemannian metric $g$ on $M$ and assume that $\JJJ$ is compatible with $g$ and parallel
with respect to its Levi-Civit\`a connection. 
Let $\o:=\o_g=\i/2\pi\cdot\gij(z)\dzidzjb$ denote its corresponding \K form, 
a closed positive $(1,1)$-form on $(M,\h{\rm J})$.
Let $H_g$ denote the Hodge projection operator from
the space of closed forms onto the kernel of $\D$. Let $V=\O^n([M])=\intm\on$.
Denote by $\HO$ the space of \K forms representing $\O$.

Let 
$PSH(M,\JJJ,\o)\sseq\Loneloc(M)$ denote the set of $\o$-plurisubharmonic functions.
Define the space of smooth strictly $\o$-plurisubharmonic 
functions (\K potentials)
$$
\calH_\o=\{\vp\in\CinfM\,:\, \ovp:=\o+\i\ddbar\vp >0\}.
$$
Let $\Ric\omega=-\i/2\pi\cdot\ddbar\log\det(\gij)$ denote the
Ricci form of $\o$. It is well-defined globally and represents
the first Chern class $c_1:=c_1(T^{1,0}M,\h{\rm J})\in H^2(M,\ZZ)$.
Alternatively it may be viewed as minus the curvature form of the canonical line bundle 
$K_M$, the top exterior product of the holomorphic cotangent bundle $T^{1,0\,\star}M$.
One calls $(M,\JJJ)$ Fano when $c_1>0$.
One calls $\o$ \KE if $\Ric\o=a\o$ for some real $a$. 
Let $f_{\o}\in\CinfM$ denote 
the unique function satisfying $\i\ddbar f_{\o}=\Rico-\o$ 
and $\Vm\intm e^{f_{\o}}\on=1$.

Let $\Aut(M,\JJJ)_0$ denote the identity component of the complex Lie group $\Aut(M,\JJJ)$ of automorphisms (biholomorphisms)
of $(M,\JJJ)$
and denote by  $\aut(M,\JJJ)$ its Lie algebra of infinitesimal automorphisms composed
of real vector fields $X$ satisfying $\L_X\JJJ=0$. 
Denote by $\H_\o(G)\sseq\Ho$ the subspace of $G$-invariant potentials.

For $\vp\in PSH(M,\JJJ,\o)$ define the multiplier ideal sheaf associated to $\vp$ to be the
sheaf $\calI(\vp)$ defined for each open set $U\sseq M$ by local sections
$$
\calI(\vp)(U)=\{h\in \calO_M(U): |h|^2 e^{-\vp}\in\Loneloc(M) \}.\eqn\IdealSheafSectionsEq
$$
Such sheaves are coherent \ref{D, p. 73; N}. 
Such a sheaf is called proper if it is neither zero nor the structure sheaf $\calO_M$.

\FDef{\ref{N, Definition 2.4}
\TaggDef{NadelSheafDef}
A proper multiplier ideal sheaf $\calI(\vp)$ will be called a Nadel sheaf whenever there 
exists $\eps>0$ such that $(1+\eps)\vp\in PSH(M,\JJJ,\o)$.
}
\noindent
According to \Ref{NadelVanishingThm} such sheaves satisfy \UsualNadelVanishingEq.
Define the complex singularity exponent $c_\o(\vp)$ of a function $\vp\in PSH(M,\JJJ,\o)$ by
$$
c_{\o}(\vp)=\sup\{\,\gamma\,:\, \intm e^{-\gamma\vp}\on<\infty\}.\eqn\ComplexSingularityExponentEq
$$
Note that $\intm e^{-\gamma\vp}\on=\infty$ implies that $\intm e^{-(\gamma+\eps)\vp}\on=\infty$ for any $\eps>0$.

Denote by $\lfloor x\rfloor$ the largest integer not larger than $x$.

\subsection{Proof of the main result}
\TaggSection{MainSection}%
The proof is split into steps similarly to the setting of the 
continuity method. First, we obtain an upper bound for
$\V\intm-\vpt\ovptn$ in terms of $\V\intm\vpt\on$.
Second, we show that if the complex singularity exponents of the functions
$\vpt$ are uniformly bigger than $n/(n+1)$ then one has a uniform upper bound on
$\V\intm\vpt\on$, and hence on $\sup\vpt$. Third, we show that $-\inf\vpt$ is
uniformly bounded from above in terms of $\V\intm-\vpt\ovptn$.
Fourth, we construct the Nadel multiplier ideal sheaf.

We turn to the proof. 
We assume throughout that $n>1$ (the remaining case will be handled separately at the end).
First we recall some necessary facts and estimates.

Consider the Ricci flow on a Fano manifold
$$
\eqa{
\frac{\pa \o(t)}{\pa t} & =-\Ric\o(t)+\o(t),\quad t\in\RR_+,\cr
\o(0) & =\o,
}\eqn\KRFEq$$
and a corresponding flow equation on the space of \K potentials $\Ho$,
$$
\ovptn=\on e^{\fo-\vpt+\dot\vpt},\q\,\vp(0)=c_0,\quad t\in\RR_+.\eqn\FlowPotentialMAEq
$$
This flow exists for all $t>0$ \ref{Ca}. Here 
$c_0$ is a constant uniquely determined by $\o$ \ref{PSS, (2.10)}.
This choice is necessary in order to have the second estimate of \Ref{RFEstimatesThm} (i) below
\ref{CT,\S\S10.1; PSS,\S2}. 

Let 
$$
||\psi||_{L^p(M,\o)} = \Big(\V\intm\psi^p \on\Big)^{\frac1p},
$$
and let
$$
\eqa{
||\psi||^2_{W^{1,2}(M,\o)} & =||\N\psi||^2_{L^2(M,\o)}+||\psi||^2_{L^2(M,\o)}
\cr & =
\V\intm n\i\del\psi\w\dbar\psi\w\o^{n-1}+\V\intm\psi^2\on. 
}
$$

We will make essential use of the following estimates of Perelman, Ye and Zhang.%

\FThm{
\TaggThm{RFEstimatesThm}
Let $(M,\JJJ)$ be a Fano manifold of complex dimension $n>1$ and let $\vpt$ satisfy \FlowPotentialMAEq.
There exist $C_1, C_2>0$ depending only on $(M,\JJJ,\o)$ such that:\newline
(i) \ref{ST,TZ} One has 
$$
||\fovpt||_{\Linf(M)}\le C_1,\q ||\dot\vpt||_{\Linf(M)}\le C_1, \q\all\, t>0.\eqn\PerelmanIneq
$$
\newline
(ii) \ref{Ye,Z} For all $\psi\in W^{1,2}(M,\o)$ one has
$$
||\psi||_{L^{\frac{2n}{n-1}}(M,\ovpt)}\le C_2 ||\psi||_{W^{1,2}(M,\ovpt)}, \q\all\, t>0.\eqn\SobolevIneq
$$
}

\vglue-0.19cm
Following Aubin \ref{A2}, define functionals on $\Hc\times\Hc$ by
$$\eqano{
I(\o,\ovp) & =
\Vm\int_M\i\del\vp\w\dbar\vp\w\sum_{l=0}^{n-1}\o^{n-1-l}\w\ovp^{l} =
\Vm\intM\vp(\on-\ovpn),&\eqnno\Ieq\cr
J(\o,\ovp) & =\frac{\Vm}{n+1}\int_M\i\del\vp\w\dbar\vp\w\sum_{l=0}^{n-1}(n-l)\o^{n-l-1}\w\ovp^{l}.
&\eqnno\Jeq
}$$
Following Ding  \ref{Di}, define a functional on $\Hc\times\Hc$ by
$$
F(\o,\ovp)
=
-(I-J)(\o,\ovp)-\V\intm\vp\ovpn
-\log\V\intm e^{f_\o-\vp}\on.\eqn\FFunctionalEq
$$
It is exact, that is to say it satisfies the cocycle condition
$F(\o_1,\o_2)+F(\o_2,\o_3)=F(\o_1,\o_3)$.
Its critical points are precisely \KE metrics.
The following monotonicity result is well-known (see, e.g., \ref{CT, Lemma 3.7}).
\FLem{
\TaggLemm{FMonotonicityLemm}
The functional $F$ is monotonically decreasing along the flow \KRFEq.
}

\noindent{\it (i) First step.}
As a consequence of \Ref{FMonotonicityLemm} and \PerelmanIneq\ we have
$$\eqa{
0\ge F(\o,\ovpt)
& =
-(I-J)(\o,\ovpt)-\V\intm\vpt\ovptn
-\log\V\intm e^{-\dot\vpt}\ovptn.
\cr
& \ge
-(I-J)(\o,\ovpt)-\V\intm\vpt\ovptn-C_1.
}\eqn\FDescreasingEq
$$

From \Ieq-\Jeq\ it follows that $\frac1{n+1}I\le J$. We then have
$$\eqa{
\V\intm-\vpt\ovptn
& \le (I-J)(\o,\ovpt)+C_1
\cr & 
\le \frac n{n+1}I(\o,\ovpt)+C_1=\frac n{n+1}\V\intm\vpt(\on-\ovptn)+C_1.
}$$
Hence,
$$
\V\intm-\vpt\ovptn
\le
\frac nV\intm\vpt\on+(n+1)C_1.\eqn\FirstStepIneq
$$
This completes the first step of the proof.

\noindent {\it (ii) Second step.}
Assume $\gamma>0$ is such that 
$$
\V\intm e^{-\gamma(\vpt-\V\intm\vpt\on)}\on\le C,
$$
with $C$ independent of $t$.
Using the flow equation, rewrite this as
$$
\V\intm e^{(1-\gamma)\vpt+\gamma\V\intm\vpt\on-\dot\vpt-f_\o}\ovptn\le C.
$$
Jensen's inequality gives,
$$
\V\intm \Big((1-\gamma)\vpt+\gamma\V\intm\vpt\on-\dot\vpt-f_\o\Big)\ovptn\le \log C.
$$
Using \PerelmanIneq\ and \FirstStepIneq\ we obtain
$$
\gamma\V\intm\vpt\on\le (1-\gamma)\V\intm -\vpt\ovptn+C'
\le n(1-\gamma)\V\intm \vpt\o^n+C''.
$$
Whenever $\gamma\in(\frac n{n+1},1)$ this yields an a priori estimate
on $\V\intm\vpt\on$. 

Under this assumption this also implies an a priori upper bound on $\sup\vpt$. Indeed, 
let $G_\o$ be a Green function for $-\D=-\D_\o$ 
satisfying $\intM G_\o(\cdot,y)\on(y)=0$. Set $A_\o=-\inf_{M\times M} G_\o$.
Recall the sub-mean value property of \opsh functions:
$$
\vp(p)\le \V\intm\vp\on+nA_\o,\q \all\,p\in M.\eqn\SupBound
$$
In fact, by assumption $-n<\D\vp$. The Green formula gives,
$$
\eqano{
\vp(p)-\VintM \vp\on & = 
-\VintM G_\o(p,y)\D\vp(y)\on(y)
\cr & =
\VintM (G_\o(p,y)+A_\o)(-\D\vp(y))\on(y)
\le nA_\o,\cr
}
$$
by the normalization of $G_\o$.
This completes the second step.

\noindent {\it (iii) Third step.}
This step follows from an argument used by Tian \ref{T1} for the continuity method. It adapts to
our setting thanks to \Ref{RFEstimatesThm}.

Put $\eta=\sup_M\vpt-\vpt+1$ and let $p>0$.
The first part of the argument involves reducing the $\Linf(M)$ estimate for $\eta$ to
an $L^2(M,\ovpt)$ estimate. First, 
$$\eqano{
\intM\eta^p\ovptn
& \ge
\intM\eta^p(\ovptn-\ovpt^{n-1}\w\o)
=
-\intM\eta^p\i\ddbar\eta\w\ovpt^{n-1}
\cr & =
\intM\i\del(\eta^p)\w\dbar\eta\w\ovpt^{n-1}
\cr & =
\frac{4p}{(p+1)^2}\intM\i\del\eta^{\frac{p+1}2}\w\dbar\eta^{\frac{p+1}2}\w\ovpt^{n-1}.&\eqnno\FirstMoserIneq
}$$
Combined with \SobolevIneq\ this gives
$$
\frac1{C_2^2}
||\eta||_{L^{\frac{(p+1)n}{n-1}}(M,\ovpt)}^{p+1}
\le
\frac{n(p+1)^2}{4p}
||\eta||_{L^p(M,\ovpt)}^p
+
||\eta||_{L^{p+1}(M,\ovpt)}^{p+1}.\eqn\WeightedIneq
$$
Following Tian, a Moser iteration argument \ref{T1, p. 235-236}
now allows us to conclude that
there exists a constant $C$ depending only on $C_2$ and $(M,\JJJ,\o)$ such that
$$
\sup\eta\le C ||\eta||_{L^2(M,\ovpt)}.
$$

The second part of the argument requires a uniform Poincar\'e inequality in order
to bound the $L^2(M,\ovpt)$ norm of $\eta$ in terms of its $L^1(M,\ovpt)$ norm.
Recall the following weighted Poincar\'e inequality \ref{F,\S2.4} (see also \ref{TZ}). 
\FLem{%
\TaggLemm{WeightedPoincareLemm}%
\ Let $(M,\JJJ)$ be a Fano manifold. Then for any function $\psi\in W^{1,2}(M,\o)$,
$$
\V\intm\big(\psi-\V\intm\psi e^{\fo}\on\big)^2e^{\fo}\on
\le
\V\intm n\i\del\psi\w\dbar\psi\w e^{\fo}\o^{n-1}.
$$
}

As a corollary we have, thanks to \PerelmanIneq, a uniform Poincar\'e-type inequality
along the flow:
$$
e^{-2C_1}||\eta||^2_{L^2(M,\ovpt)}
-
e^{C_1}
||\eta||^2_{L^1(M,\ovpt)}
\le
||\N\eta||^2_{L^2(M,\ovpt)}.\eqn\UniformPoincareIneq
$$

Therefore, applying \FirstMoserIneq, now with $p=1$, combined with \UniformPoincareIneq, we obtain
$$
e^{-2C_1}||\eta||^2_{L^2(M,\ovpt)}
-
e^{C_1}
||\eta||^2_{L^1(M,\ovpt)}
\le
n||\eta||_{L^1(M,\ovpt)},
$$
which completes the second part of the argument.

Finally,
$$
||\eta||_{L^1(M,\ovpt)}
=
1+\sup\vpt+\V\intm-\vpt\ovptn,\eqn\InfBound
$$
and this is uniformly bounded in terms of $\V\intm\vpt\on$ thanks to \FirstStepIneq\ and
\SupBound. This completes the third step.

\noindent
{\it (iv) Fourth step.} Assume that $(M,\JJJ)$ does not admit a \KE metric. The
relevant theory of complex \MA equations due to Aubin and Yau \ref{A1,Y}
now implies that 
$\{||\vpt||_{\Linf(M)}\}_{t\in[0,\infty)}$ is unbounded. If not, one would have uniform
higher-order estimates on $\vpt$ and properties of Mabuchi's K-energy \ref{M} (equivalent in many 
ways to $F$)
then show that a subsequence will converge to a smooth \KE metric (see \ref{PSS,\S2}).
Combining the first three steps 
this implies that
for each $m\in\NN$ we may find an unbounded increasing subsequence of times $\{t_{j(m)}\}_{j(m)\ge1}\sseq[0,\infty)$ for which 
$$
\lim_{j\ra\infty}\intM  e^{-(\frac n{n+1}+\frac1m)\big(\vp_{t_{j(m)}}-\V\intm\vp_{t_{j(m)}}\on\big)}\on=\infty.
$$
Hence, by the diagonal argument, there exists a subsequence of potentials 
$\{\vp_{t_j}\}_{j\ge1}$
for which (one may equivalently work throughout with $\sup\vp_{t_j}$ instead of $\V\intm\vp_{t_j}\on$)
$$
\lim_{j\ra\infty}\intM  e^{-\gamma \big(\vp_{t_j}-\V\intm\vp_{t_j}\on\big)}\on=\infty, \q\all\,\gamma\in(n/(n+1),\infty).\eqn
$$
One may further extract an unbounded increasing sub-subsequence of times such that
$\{\vp_{t_{j_k}}-\Vm\intm\vp_{t_{j_k}}\on\}_{k\ge1}$ converges in the $L^1(M,\o)$-topology
to a limit $\vp_\infty\in PSH(M,\JJJ,\o)$ \ref{DK, p. 549-550}.
The Demailly-Koll\'ar lower semi-continuity of complex singularity exponents then implies
$$
c_\o(\vp_\infty)\le \liminf_{k\ra\infty} 
c_\o(\h{$\vp_{t_{j_k}}-\Vm\intm\vp_{t_{j_k}}\on$})
\le \frac n{n+1}.\eqn
$$
Equivalently,
$$
||e^{-\gamma\vp_\infty}||_{L^1(M,\o)}=\infty,\q\all\,\gamma\in(n/(n+1),\infty),
$$ 
and the multiplier ideal sheaf $\calI(\gamma\vp_\infty)$
defined for each open set $U\sseq M$ by local sections
$\calI(\gamma\vp_\infty)(U)=\{h\in \calO_M(U): |h|^2 e^{-\gamma\vp_\infty}\in\Loneloc(M) \}$ is not 
$\calO_M$.
It is also not zero since $\vp_\infty$ is not identically $-\infty$ as 
its average is zero.

Finally, we recall a version of Nadel's vanishing theorem.

\FThm{
\TaggThm{NadelVanishingThm}
\ref{N}
Let $(M,\JJJ,\o)$ be a Hodge manifold and $(L,h)$ an ample holomorphic line bundle over $M$ equipped
with a smooth Hermitian metric with positive curvature form $\Psi_h$ given locally 
by $-\i/2\pi\cdot\ddbar\log h$. 
Assume that $(1+\eps)\vp\in PSH(M,\JJJ,\Psi_h)$ for some $\eps>0$.
Then
$$
H^r(M,\calI(\vp)\otimes K_M\otimes L)=0, \q \all\, r\ge 1.\eqn\NadelVanishingEq
$$
}

\vglue-0.5cm
The cohomology vanishing statement \NadelVanishingFanoEq\
for the sheaf $\calI(\gamma\vp_\infty)$ just constructed
is a consequence of \Ref{NadelVanishingThm}
with $L=K_M^{-\lfloor\gamma\rfloor-1}$  since  
$$
(\lfloor\gamma\rfloor+1) \o+\gamma\i\ddbar\vp_\infty\ge
(\lfloor\gamma\rfloor+1-\gamma)\o.
$$
This concludes the proof of \Ref{MainThm} for $n>1$.

To treat the remaining case $n=1$ ($M=S^2$), 
replace the third step with the following argument. First, one 
has a uniform lower bound for the scalar
curvature along the flow (see, e.g., \ref{C}), i.e., a lower bound on the Ricci curvature
when $n=1$. Next,
by Perelman \ref{ST} a uniform diameter bound holds along the flow. Therefore the quantity
$\diam(M,g_{\ovpt})^2\cdot\inf\{\,\Ric\ovpt(v,v)\,:\, {||v||_{g_{\ovpt}}=1}\}$ is uniformly bounded from below.
Applying the Bando-Mabuchi Green's function estimate \ref{BM, Theorem 3.2} now implies that
$A_{\ovpt}$ (see \SupBound\ for notation) is uniformly bounded from above (here we invoked Perelman's diameter bound again). 
Since $\D_{\ovpt}\vpt<1$, Green's formula gives
a uniform bound for $-\inf\vpt$ in terms of $\V\intm-\vpt\ovpt$. Now the first, second and fourth steps
apply without change to give the desired result. \done

\vglue-0.3cm
\subsection{Remarks and further study}
\TaggSection{RemarksSection}
We end with some 
remarks.

We saw that the uniform Sobolev inequality of Ye and Zhang can be used instead of
the \Kolodziej\ Theorem in order to obtain the $\Linf$ estimate for \Ref{ConvergenceCor}. 
Observe that this also applies to the proof of the main theorem in \ref{TZ}, at least
in the case of no holomorphic vector fields. 
Indeed, by \SupBound, \InfBound\ and \FirstStepIneq, $\sup\vpt$ and $-\inf_{\vpt}$ are bounded 
in terms of $\V\int\vpt\on$. By   
\FlowPotentialMAEq\ and \PerelmanIneq\ one has
$\V\int\vpt\on\le I(\o,\ovpt)+C_1$. Hence, if a \KE metric exists the properness of $F$ in the
sense of Tian \ref{T3} 
implies a uniform $\Linf$ bound on $\vpt$ (compare \ref{PSS,\S3}). The same applies
also to Pali's theorem \ref{P}.

The statement of \Ref{NadelThm} can be refined to hold for all
$\gamma\in(t_0n/(n+1),\infty)$ with $t_0=t_0(\o)\le 1$ defined to be the first time
for which $\{||\vpt||_{\Linf(M)}\}_{t\in[0,t_0)}$ is unbounded,
with $\vpt$ a solution of  \AubinPathEq\
\ref{T1, p. 234} (see also \ref{N, p. 582}).
One has the inequality $t_0\le \sup
\{\,b\,:\,\Ric\o\ge b\o, \o\in\Hc,b\in\RR\}$, the right hand side a holomorphic 
invariant studied by Tian; on some Fano manifolds it is smaller than 1 \ref{T2}.
We do not know whether for these manifolds the exponent in 
\Ref{MainThm} can be lowered as well. Perhaps the difference here is related 
to the fact that the flow always exists for all
time unlike the continuity path. 
Yet, as far as the usefulness of the sheaves is concerned, this does not seem to be crucial since
they all satisfy the same vanishing conditions \UsualNadelVanishingEq\ for exponents smaller than
1.
It is worth mentioning that this invariant is smaller than 1 only when the functional $F$ 
is not bounded from below \ref{BM,DT}. It would be very  
interesting to know what can be said 
regarding the converse (compare \ref{R1,\S1}). We are therefore
led to pose the following problem:

\FProb{On a Fano manifold, determine whether the lower boundedness of the functional $F$ (equivalently, 
of Mabuchi's K-energy)
on $\Hc$ is equivalent to $\sup
\{\,b\,:\,\Ric\o\ge b\o, \o\in\Hc,b\in\RR\}=1$.}

The \KR flow has been widely studied for manifolds whose first Chern class is definite or zero (see, e.g., 
\ref{Ch2}).
This flow is simply the
dynamical system induced by integrating minus the Ricci potential vector field $-f$ on the space
of \K forms $\HO$.
A vector field $\psi$ on the space of \K forms is an assignment 
$\o\mapsto \psi_\o\in \CinfM/\RR$.
The vector field $f$ is the assignment $\o\mapsto f_\o$ with $\fo$ defined by 
$\Ric\o-\mu\o=\i\ddbar\fo$, $\,\mu\in\RR$. 

Motivated by this observation one is naturally led to extend the definition of the \KR flow to an arbitrary \K
manifold, simply by defining the flow lines to be integral curves of minus the Ricci potential vector field $-f$
on $\HO$, with $\O$ an arbitrary \K class. 
Recall that the Ricci potential is defined in general by 
$\Ric\o-H_\o\Ric\o=\i\ddbar\fo$. The resulting flow equation can also be written as
$$
\eqa{
\frac{\pa \o(t)}{\pa t} & =-\Ric\o(t)+H_t\Ric\o(t),\quad t\in\RR_+,\cr
\o(0) & =\o,
}\eqn\KRFEq$$
for each $t$ for which a solution exists in $\HO$. This flow, introduced by Guan, 
is part of the folklore in the field although it has not been much studied.\note{%
\cmreight It seems that Guan first considered this flow in unpublished work in
the 90's (see references to
\reffoot{G1}).
After posting the first version of this note I became aware, thanks to
G. Sz\'ekelyhidi, of a recent preprint \reffoot{G2} posted by Guan on his webpage in which
this flow is studied. For a different but related flow see \reffoot{S}.
}

Several authors have raised the question whether Nadel's construction can be extended to the study
of constant scalar curvature \K metrics. We believe that Nadel-type obstruction sheaves should arise 
from this dynamical system (as well as from its discretization \ref{R2,R3}; in these two references
a ``discrete" analogue of \Ref{PSSThm} was shown to hold) in the absence of fixed points. 

As we saw, it is important to make a choice of how to induce a flow on $\Ho$ from that on $\HO$,
and different normalizations give rise to different sheaves.
The flow equation \FlowPotentialMAEq\ corresponds to restricting the evolution 
to a certain codimension one submanifold of $\Ho$. For a general \K class one may
define an operator on the space of \K forms $\HO$, identified as an open subset of $\CinfM/\RR$, by
$$
h:\vp\in\HO\mapsto h(\vp)\in\CinfM/\RR,\eqn
$$
with 
$$
H_{\ovp}\Ric\ovp-H_\o\Ric\o =\i\ddbar h(\vp).\eqn
$$
By choosing such an appropriate submanifold,
analogously to \FlowPotentialMAEq, one may consider the flow
$$
\o_{\vpt}^n=\on e^{f_\o-h(\vpt)+\dot\vpt}\on.\eqn\GeneralContinuityPathEq
$$
The difficulty lies in the fact that this operator is in general no longer a multiple of the identity.

\def\smallblackbox{\vrule height.6ex width .5ex depth -.1ex}
\def\boxseparation{\hfil\smallblackbox$\q$\smallblackbox$\q$\smallblackbox\hfil}
%\vbox{\bigskip
\boxseparation
\bigskip

This note was written during a visit to the Technion and I thank that institution,
and especially M. Cwikel and Y. Pinchover, for the hospitality and partial financial support.
I am indebted to my teacher, Gang Tian, for his advice and 
warm encouragement.
I am grateful to Yum-Tong Siu, from whose class I first learned on some of Nadel's work.
I thank D. Kim for enjoyable discussions on multiplier ideal sheaves, 
V. Tosatti and a referee for valuable comments on this note, and Q. Zhang for a useful discussion on 
his article \ref{Z}.
This material is based upon work supported under a National Science 
Foundation Graduate Research Fellowship.
%}

\frenchspacing

\bigskip\bigskip
\noindent{\bf Bibliography}
\bigskip
\def\ref#1{\Taggf{#1}\item{ {\bf[}{\sans #1}{\bf]} } }

\def\sm{\vglue1.8pt}

\ref{A1} Thierry Aubin, \'{E}quations du type {M}onge-{A}mp\`ere sur les vari\'et\'es
k\"ahl\'eriennes compactes, {\sl Bulletin des Sciences Math\'ematiques} 
{\bf 102} (1978), 63--95.
\sm
\ref{A2} \opcit, R\'eduction du cas positif de l'\'equation de {M}onge-{A}mp\`ere 
sur les vari\'et\'es k\"ahl\'eriennes compactes \`a la d\'emonstration d'une in\'egalit\'e, 
{\sl Journal of Functional Analysis} {\bf 57} (1984), 143--153.
\sm
\ref{BM} Shigetoshi Bando, Toshiki Mabuchi, Uniqueness of K\"ahler-Einstein metrics
modulo connected group actions, in {\it Algebraic Geometry,
Sendai, 1985} (T. Oda, Ed.), Advanced Studies in Pure Mathematics {\bf 10},
Kinokuniya, 1987, 11--40.
\sm
\ref{Ca} Huai-Dong Cao, Deformations of K\"ahler metrics to K\"ahler-Einstein metrics on compact
\K manifolds, {\sl Inventiones Mathematicae} {\bf 81} (1985), 359--372.
\sm
\ref{CLT} Xiu-Xiong Chen, Peng Lu, Gang Tian, A note on uniformization of Riemann surfaces
by Ricci flow, {\sl Proceedings of the American Mathematical Society} {\bf 134} (2006), 3391--3393.
\sm
\ref{CT} Xiu-Xiong Chen, Gang Tian, Ricci flow on K\"ahler-Einstein surfaces,
{\sl Inventiones Mathematicae} {\bf 147} (2002), 487--544.
\sm
\ref{C} Xiu-Xiong Chen, On K\"ahler manifolds with positive orthogonal bisectional curvature,
{\sl Advances in Mathematics} {\bf 215} (2007), 427--445.
\sm
\ref{Ch1} Bennett Chow, The Ricci flow on the 2-sphere, {\sl Journal of Differential Geometry}
{\bf 33} (1991), 325--334.
\sm
\ref{Ch2} Bennett Chow et al., The Ricci flow: Techniques and applications. Part I: Geometric
aspects, American Mathematical Society, 2007.
\sm
\ref{D} Jean-Pierre Demailly, On the Ohsawa-Takegoshi-Manivel $L^2$ extension theorem,
in {\it Complex analysis and Geometry}  (P. Dolbeault et al., Eds.), 
Progress in Mathematics {\bf 188}, Birkh\"auser, 2000, 47--82.
\sm
\ref{DK} Jean-Pierre Demailly, J\'anos Koll\'ar, Semi-continuity of complex singularity exponents
and K\"ahler-Einstein metrics on Fano orbifolds, {\sl Annales Scientifiques de l'\'Ecole
Normale sup\'erieure} {\bf 34} (2001), 525--556.
\sm
\ref{Di} Wei-Yue Ding, Remarks on the existence problem of positive
{K}\"ahler-{E}instein metrics, {\sl Mathematische Annalen} {\bf 282} (1988), 463--471.
\sm
\ref{DT} Wei-Yue Ding, Gang Tian, The generalized Moser-Trudinger inequality, 
in {\it Nonlinear Analysis and Microlocal Analysis: Proceedings of the International
Conference at Nankai Institute of Mathematics} (K.-C. Chang et al., Eds.), 
World Scientific, 1992, 57--70. ISBN 9810209134.
\sm
\ref{F} Akito Futaki, K\"ahler-{E}instein metrics and integral invariants,
Lecture Notes in Mathematics {\bf 1314}, Springer, 1988.
\sm
\ref{G1} Daniel Z.-D. Guan, Quasi-Einstein metrics, {\sl International Journal of Mathematics}
{\bf 6} (1995), 371--379.
\sm
\ref{G2} \opcit, Extremal-solitons and $C^\infty$ convergence of the modified Calabi flow on certain
$CP^1$ bundles, preprint, 2006.
\sm
\ref{H1} Richard S. Hamilton, Three-manifolds with positive Ricci curvature,
{\sl Journal of Differential Geometry} {\bf 17} (1982), 255--306.
\sm
\ref{H2} \opcit, The Ricci flow on surfaces, in {\it Mathematics and general relativity} (J. I. Isenberg, Ed.),
Contemporary Mathematics {\bf 71}, American Mathematical Society, 1988, 237--262.
\sm
\ref{K} Joseph J. Kohn, Subellipticity of the $\bar \partial$-Neumann problem on pseudo-convex domains: sufficient conditions, Acta Mathematica {\bf 142} (1979), 79--122.
\sm
\ref{Ko} S\polishl awomir Ko\polishl odziej, The complex Monge-Amp\`ere equation, 
{\sl Acta Mathematica} {\bf 180} (1998), 69--117.
\sm
\ref{M} Toshiki Mabuchi, K-energy maps integrating Futaki invariants,
{\sl T\^ohoku Mathematical Journal} {\bf 38} (1986), 575--593.
\sm
\ref{N} Alan M. Nadel, Multiplier ideal sheaves and K\"ahler-Einstein metrics of positive scalar
curvature, {\sl Annals of Mathematics} {\bf 132} (1990), 549--596.
\sm
\ref{P} Nefton Pali, Characterization of Einstein-Fano manifolds via the K\"ahler-Ricci flow,
preprint, arxiv:math.DG/0607581v2.
\sm
\ref{PSS} Duong H. Phong, Nata\v sa \v Se\v sum, Jacob Sturm,
Multiplier ideal sheaves and the K\"ahler-Ricci flow, preprint,
arxiv:math.DG/0611794v2. To appear in {\sl Communications in Analysis and Geometry}.
\sm
\ref{R1} Yanir A. Rubinstein, On energy functionals, K\"ahler-Einstein metrics,
and the Moser-Trudinger-Onofri neighborhood, preprint, arxiv:math.DG/0612440. To appear in
{\sl Journal of Functional Analysis}.
\sm
\ref{R2} \opcit, The Ricci iteration and its applications, 
Comptes Rendus de l'Acad\'emie des Sciences Paris, Ser. I {\bf 345} (2007), 445--448.
\sm
\ref{R3} \opcit, Some discretizations of geometric evolution equations and 
the Ricci iteration on the space of K\"ahler metrics, I, preprint, arxiv:0709.0990 [math.DG].
\sm
\ref{ST} Nata\v sa \v Se\v sum, Gang Tian, Bounding scalar curvature and diameter 
along the K\"ahler-Ricci flow (after Perelman) and some applications, preprint.
\sm
\ref{S} Santiago R. Simanca, Heat flows for extremal K\"ahler metrics, 
{\sl Annali della Scuola Normale Superiore di Pisa} {\bf 4} (2005), 187--217.
\sm
\ref{Si} Yum-Tong Siu, The existence of K\"ahler-Einstein metrics on manifolds with 
positive anticanonical line bundle and a suitable finite symmetry group, 
{\sl Annals of Mathematics} {\bf 127} (1988), 585--627.
\sm
\ref{So} Jian Song, The $\alpha$-invariant on toric Fano manifolds,
{\sl American Journal of Mathematics} {\bf 127} (2005), 1247--1259.
\sm
\ref{T1} Gang Tian, On K\"ahler-Einstein metrics on certain K\"ahler manifolds with $C_1(M)>0$,
{\sl Inventiones Mathematicae} {\bf 89} (1987), {225--246}.
\sm
\ref{T2} \opcit, On stability of the tangent bundles of Fano varieties,
{\sl International Journal of Mathematics} {\bf 3} (1992), 401--413.
\sm
\ref{T3} \opcit,
K\"ahler-{E}instein metrics with positive scalar curvature, 
{\sl Inventiones Mathematicae} {\bf 130} (1997), {1--37}.

\sm
\ref{TY} Gang Tian, Shing-Tung Yau, K\"ahler-Einstein metrics on complex surfaces
with $C_1>0$, {\sl Communications in Mathematical Physics} {\bf 112} (1987), 175--203.
\sm
\ref{TZ} Gang Tian, Xiao-Hua Zhu, Convergence of K\"ahler-Ricci flow, 
{\sl Journal of the American Mathematical Society} {\bf 20} (2007), 675--699.
\sm
\ref{Y} Shing-Tung Yau, On the Ricci curvature of a compact K\"ahler
manifold and the Complex Monge-Amp\`ere equation, I, {\sl Communications in Pure
and Applied Mathematics} {\bf 31} (1978), 339--411.
\sm
\ref{Ye} Rugang Ye, The logarithmic Sobolev inequality along the Ricci flow, preprint,
arxiv:0707.2424v4 [math.DG].
\sm
\ref{Z} Qi S. Zhang, A uniform Sobolev inequality under Ricci flow,
{\sl International Mathematics Research Notices} (2007), Article ID rnm056. Erratum
to: A uniform Sobolev inequality under Ricci flow (2007), Article
ID rnm096.

\end